# Introduction of Over/Under/Off Masses

**Florentin Smarandache** [1,*]

[1] Math and Science Division, University of New Mexico, Gallup, NM 87301, U.S.A.

**\*** Correspondence: smarand@unm.edu

**Abstract:** In this paper, for the first time, we extend the Over/Under/Off Set/Logic/Probability used in uncertain theories (such as: fuzzy, neutrosophic and extensions) to the Over/Under/Off Mass that could be used in information fusion. The approach is exemplified in three scenarios: (1) wildfire evacuation and resource allocation with satellite, IoT, and social media; (2) coverage gaps where indeterminacy must be managed; and (3) security monitoring where contradictory or erroneous reports are discounted.

## 1. Introduction

By "uncertain" set/logic we mean all types of fuzzy and fuzzy-extensions (intuitionistic fuzzy, neutrosophic, refined neutrosophic, spherical fuzzy, plithogenic, etc.) of set and logic. The uncertain Set was extended by Smarandache in 2007 to uncertain OverSet (when some component is > 1), since he observed that, for example, an employee working overtime deserves a degree of membership > 1, with respect to an employee that only works regular full-time and whose degree of membership = 1; and to uncertain UnderSet (when some neutrosophic component is < 0), since, for example, an employee making more damage than benefit to his company deserves a degree of membership < 0, with respect to an employee that produces benefit to the company and has the degree of membership > 0; and to uncertain OffSet (when some neutrosophic components are off the interval [0, 1], i.e. some neutrosophic component > 1 and some neutrosophic component < 0). As a particular case it was considered the Neutrosophic Over-/Under-/Off-Set / Logic / Measure / Probability / Statistics. More generally, one has the uncertain Set/Logic/Measure/Probability/Statistics etc. were extended to respectively uncertain Over-/Under-/Off- Set / Logic / Measure / Probability / Statistics etc.

## 2. Example of OverMembership, UnderMembership, OffMembership [1]

In a company a full-time employer works 40 hours per week. Let's consider the last week period. Helen worked part-time, only 30 hours, and the other 10 hours she was absent without payment; hence, her membership degree was 30/40 = 0.75 < 1. John worked full-time, 40 hours, so he had the membership degree 40/40 = 1, with respect to this company. But George worked overtime 5 hours, so his membership degree was (40+5)/40 = 45/40 = 1.125 > 1.

Thus, we need to make distinction between employees who work overtime, and those who work full-time or part-time. That's why we need to associate a degree of membership strictly greater than 1 to the overtime workers. Now, another employee, Jane, was absent without pay for the whole week, so her degree of membership was 0/40 = 0. Yet, Richard, who was also hired as a full-time, not only didn't come to work last week at all (0 worked hours), but he produced, by accidentally starting a devastating fire, much damage to the company, which was estimated at a value half of his salary (i.e. as he would have gotten for working 20 hours that week).

Therefore, his membership degree has to be less that Jane's (since Jane produced no damage). Whence, Richard's degree of membership, with respect to this company, was - 20/40 = - 0.50 < 0. Consequently, we need to make distinction between employees who produce damage, and those who produce profit, or produce neither damage no profit to the company. Therefore, the membership degrees > 1 and < 0 can fit with some real situations, so we have to take them into consideration. Then, similarly, the Neutrosophic Logic/Measure/Probability/Statistics etc. were extended to respectively Neutrosophic Over/Under/Off-Logic, -Measure, -Probability, - Statistics etc. [2].

## 3. Frame of Discernment and Its Powerset

This section may be divided by subheadings. It should provide a concise and precise description of the experimental results, their interpretation as well as the experimental conclusions that can be drawn.

Let $\Theta = \{\theta_1, \theta_2, \ldots, \theta_n\}$, for integer $n \geq 2$, be a frame of discernment of $n$ elements, and $2^\Theta$ the powerset of $\Theta$:

[1] $\emptyset, \theta_1, \theta_2, \ldots, \theta_n \in 2^\Theta$

[2] If $A, B \in 2^\Theta$, then $A \cup B \in 2^\Theta$.

[3] No other elements belong to $2^\Theta$, except those obtained by using rules 1 or 2.

## 4. Definition of Classical Mass of Belief [3]

$$m: 2^\Theta \to [0,1]$$

$$\sum_{x \in 2^\Theta} m(x) = 1$$

## 5. Fusion of Classical Masses (Conjunction Rule)

Let $m_1, m_2 : 2^\Theta \to [0, 1]$, be two classical masses.

Then $m_{1,2} : 2^\Theta \to [0, 1]$, for any $A \in 2^\Theta$.

$$m_{1,2}(A) = \sum_{\substack{x,y \in 2^\Theta \\ x \cap y = A}} m_1(x) \cdot m_2(y)$$

In *Dempster-Shafer Theory* (DST) one applies *Dempster's Rule* to combine two distinct sources of evidence characterized by their masses [3]:

$m_{1,2}(A) = 0$, for $A = \emptyset$,

$m_{1,2}(A) = \frac{1}{1-k} \cdot \sum_{\substack{x,y \in 2^\Theta \\ x \cap y = A}} m_{1(x)} \cdot m_2(y)$, for $A \neq \emptyset$,

where $k$ is the total conflict:

$$k = \sum_{\substack{x,y \in 2^\Theta \\ x \cap y = \emptyset}} m_1(x) \cdot m_2(y)$$

## 6. Frame of Discernment and Its Hyper-Powerset

Let $\Theta = \{\theta_1, \theta_2, \ldots, \theta_n\}$, for integer $n \geq 2$, be a frame of discernment of $n$ elements, and $G^\Theta$ the hyper-powerset of $\Theta$:

1. $\emptyset, \theta_1, \theta_2, \ldots, \theta_n \in G^\Theta$.

2. If $A, B \in 2^\Theta$, then $A \cup B \in G^\Theta$.

3. If $A, B \in 2^\Theta$, then $A \cap B \in G^\Theta$.

No other elements belong to $G^\Theta$, except those obtained by using rules 1, 2, or 3.

   A. In Dezert-Smarandache Theory (DSmT) [4] we apply the PCR5 Rule of Combination for two sources of evidence:

$$m_{PCR5}(\phi) = 0$$

and for all $A \in G^\Theta \setminus \{\phi\}$

$$m_{PCR5}(A) = m_{12}(A) + \sum_{\substack{X \in G^\Theta \setminus \{A\} \\ X \cap A = \phi}} [\frac{m_1(A)^2 \cdot m_2(X)}{m_1(A) + m_2(X)} + \frac{m_2(A)^2 \cdot m_1(X)}{m_2(A) + m_1(X)}]$$

   B. where all sets involved in formulas are in canonical form and where $G^\Theta$ corresponds to classical power set $2^\Theta$ if Shafer's model is used (i.e. all elements of $\Theta$ are exclusive and exhaustive), or to a constrained hyper-power set $D^\Theta$ if any other hybrid DSm model is used instead, or to the super-power set $S^\Theta$ if the minimal refinement $\Theta^{\text{ref}}$ of $\Theta$ is used;

   C. $m_{12}(A) \equiv m_\cap(A)$ corresponds to the conjunctive consensus on A between the s = 2 sources;

   D. and where all denominators are different from zero; if a denominator is zero, that fraction is discarded.

**7. Definition of OverMass ($m_{over}$)**

$$m_{over}: 2^\Theta \to [0, \psi], \text{ where } 1 < \psi,$$

$$\sum_{x \in 2^\Theta} m(x) = 0 + \psi = \psi.$$

**8. Definition of UnderMass ($m_{under}$)**

$$m_{under}: 2^\Theta \to [\varphi, 1], \text{ where } \varphi < 0,$$

$$\sum_{x \in 2^\Theta} m(x) = \varphi + 1.$$

**9. Definition of OffMass ($m_{off}$)**

$$m_{off}: 2^\Theta \to [\varphi, \psi], \text{ where } \varphi < 0 < 1 < \psi,$$

$$\sum_{x \in 2^\Theta} m(x) = \varphi + \psi.$$

**10. Fusion of Overmasses on The Same Interval [0,Ψ]**

Let $m_{over1}, m_{over2}: 2^\Theta \to [0, \psi]$, $\psi > 1$, be two OverMasses.

Then, the *OverConjunctive Rule* is:

$m_{over1,2}: 2^\Theta \to [0, \psi]$

$$m_{over}(A) = \frac{1}{\psi} \sum_{\substack{x,y \in 2^\Theta \\ x \cap y = A}} m_{over1}(x) \cdot m_{over2}(y)$$

for any $A \in 2^\Theta$.

The multiplier $\frac{1}{\psi}$ helps to OverNormalization of all OverMasses.

**11. Example of Fusion of Overmasses on The Same Interval [0,1.1]**

A. Two workers, A and B, are applying for a manager position at a given company. They are evaluated by two company directors $D_1$ (whose evaluation is presented below as the overmass $m_{over1}$) and $D_2$ (whose, similarly, evaluation is given below as the overmass $m_{over2}$) that follow:

|            | A    | B    | A∪B  | $k_{over1,2}$ | $m_{over1}, m_{over2}: D \to [0, 1.1]$ |
|------------|------|------|------|---------------|----------------------------------------|
| $m_{over1}$ | 0.6  | 0.3  | 0.2  |               | $D = \{A, B, A \cup B\}, A \cap B = \emptyset$ |
| $m_{over2}$ | 0.5  | 0.5  | 0.1  |               |                                        |
| $m_{over12}$ | 0.46 | 0.28 | 0.02 | 0.45          | Conjunctive Rule                       |

What is the best worker to be hired for this company?

We need to fusion ($\oplus$) the two evaluations: $m_{over1} \oplus m_{over2} = m_{over1,2}$

The fusion is computed in the following way:

$m_{over\,12}(A) == [m_{over1}(A) \cdot m_{over2}(A) + m_{over1}(A) \cdot m_{over12}(A \cup B) + m_{over2}(A) \cdot m_{over1}(A \cup B)]$

$= [0.6(0.5) + 0.6(0.1) + 0.5(0.2)] = 0.46$

$m_{over\,12}(B) = 0.3(0.5) + 0.3(0.1) + 0.5(0.2) = 0.26$

$m_{over\,12}(A \cup B) = 0.2(0.1) = 0.02$

The total OverConflict:

$k_{over\,1,2} = m_{over1}(A) \cdot m_{over2}(B) + m_{over2}(A) \cdot m_{over1}(B)$

$= 0.6(0.5) + 0.5(0.3) = 0.30 + 0.15 = 0.45$

But $m_{over\,12}(A) + m_{over\,12}(B) + m_{over\,12}(A \cup B) + k_{12}$

$= 0.46 + 0.28 + 0.02 + 0.45 = 1.21$

because

$$\left[\sum_{x \in 2^\theta} m_1(x)\right] \cdot \left[\sum_{x \in 2^\theta} m_2(x)\right] =$$

$[0.6 + 0.3 + 0.2] \cdot [0.5 + 0.5 + 0.1] = [1.1] \cdot [1.1] = 1.21$

We OverNormalize by dividing by 1.1, and we get:

$$\frac{0.46}{1.1} + \frac{0.28}{1.1} + \frac{0.02}{1.1} + \frac{0.45}{1.1} = \frac{1.21}{1.1}$$

$0.418 + 0.255 + 0.018 + 0.409 = 1.1$

We use the PCR5 (Proportional Conflict Redistribution Rule #5) for the total conflict overmass $k_{over1,2} = 0.409$ to be redistributed to $A, B$, and $A \cup B$, proportionally to their masses; where $x_A, x_B, x_{A \cup B}$ represent the parts of the overmasses taking from the total conflict overmass (0.409) and added to $m_{over1,2}(A), m_{over1,2}(B)$, and $m_{over1,2}(A \cup B)$ respectively:

$$\frac{x_A}{0.418} = \frac{x_B}{0.255} = \frac{x_{A \cup B}}{0.018} = \frac{0.409}{0.418 + 0.255 + 0.018} = \frac{0.409}{0.691}$$

Hence:

$\frac{x_A}{0.418} = \frac{0.409}{0.691}$, whence $x_A = \frac{0.418(0.419)}{0.691} \simeq 0.2534616498\ldots$

$\frac{x_B}{0.255} = \frac{0.409}{0.691}$, whence $x_B = \frac{0.255(0.409)}{0.691} \simeq 0.1509334298\ldots$

$\frac{x_{A \cup B}}{0.018} = \frac{0.409}{0.691}$, whence $x_{A \cup B} = \frac{0.018(0.409)}{0.691} \simeq 0.0106541245\ldots$

whence:

$$m_{over1,2}(A) = 0.418 + 0.2534616498\ldots \simeq 0.6714616498 \simeq 0.67$$
$$m_{over1,2}(B) = 0.255 + 0.1509334298\ldots \simeq 0.4059334298 \simeq 0.40$$
$$m_{over1,2}(A \cup B) = 0.018 + 0.0106541245\ldots \simeq 0.0286541245 \simeq 0.03$$

Therefore, the company will promote worker *A* to the level of manager.

## 12. Example on different intervals [0, 1.1] and [0, 1.2]

The same example as above, [Two workers, A and B, are applying for a manager position at a given company. They are evaluated by two company directors $D_1$ (whose evaluation is presented below as the overmass $m_{over1}$) and $D_2$ (whose, similarly, evaluation is given below as the overmass $m_{over2}$) that follows, but the directors' ($D_1$ and $D_2$) evaluations are different:

|  | A | B | A ∪ B | $\phi$ |
|---|---|---|---|---|
| $m_{over1}$ | 0.7 | 0.3 | 0.1 |  |
| $m_{over2}$ | 0.4 | 0.6 | 0.2 |  |
| $m_{over1,2}$ | 0.46 | 0.30 | 0.02 | 0.54 |

$$m_{over1} : 2^\Theta \to [0, 1.1]$$
$$m_{over2} : 2^\Theta \to [0, 1.2]$$

Then

$$m_{over1} \oplus m_{over2} : 2^\Theta \to [0, 1.2]$$

$$m_{over1,2}(A) = m_{over1}(A) \cdot m_{over2}(A) + m_{over1}(A) \cdot m_{over2}(A \cup B) + m_{over2}(A \cup B) \cdot m_{over2}(A)$$
$$= 0.7(0.4) + 0.7(0.2) + 0.4(0.1) = 0.46.$$

$$m_{over1,2}(B) = m_{over1}(B) \cdot m_{over2}(B) + m_{over1}(B) \cdot m_{over2}(A \cup B) + m_{over2}(B) \cdot m_{over1}(A \cup B)$$
$$= 0.3(0.6) + 0.3(0.2) + 0.6(0.1) = 0.30.$$

$$m_{over1,2}(A \cup B) = m_{over1}(A \cup B) \cdot m_{over2}(A \cup B) = 0.1(0.2) = 0.02$$

The conflicting mass:

$$m_{over1,2}(\phi) = m_{over1}(A) \cdot m_{over2}(B) + m_{over2}(A) \cdot m_{over1}(B) = 0.7(0.6) + 0.4(0.3) = 0.54 \cdot$$

But the conflicting overmass: $m_{over\,1,2}(\phi) = 0.491$ has to be redistributed to the elements A and B involved into the conflict, proportionally with respect to their masses, according to the PCR5 Rule from Dezert-Smarandache Theory of Information Fusion.

$$\frac{X_A}{0.7} - \frac{X_B}{0.6} = \frac{0.7(0.6)}{0.7 + 0.6} = \frac{0.42}{1.3}$$

$$X_A = 0.7 \cdot \frac{0.42}{1.3} \approx 0.2261538462 \ldots$$

$$X_B = 0.6 \cdot \frac{0.42}{1.3} \approx 0.1938461538 \ldots$$

Approximate:

$$X_A \simeq 0.226, X_B \simeq 0.194.$$

Also:

$$\frac{Y_A}{0.4} - \frac{Y_B}{0.3} = \frac{0.4(0.3)}{0.4 + 0.3} = \frac{0.12}{0.7}$$

$$Y_A = 0.4 \cdot \frac{0.12}{0.7} = 0.068574286 \ldots \simeq 0.069$$

$$Y_B = 0.3 \cdot \frac{0.12}{0.7} = 0.0514285714 \ldots \simeq 0.051$$

Approximate:

$$Y_A \simeq 0.069, Y_B \simeq 0.051.$$
$$m_{\text{over1,2}}(A) = 0.46 + 0.226 + 0.069 = 0.755$$
$$m_{\text{over1,2}}(B) = 0.30 + 0.194 + 0.051 = 0.545$$

We got:

|  | A | B | A∪B | $\phi$ |
|---|---|---|---|---|
| $m_{\text{over1,2}}$ | 0.755 | 0.545 | 0.020 | 0 (zero) |

The total sum is:

0.755+0.545+0.020=1.32=1.1(1.2).

But we need to overnormalize to [0, 1.2], since

$$m_{\text{over1,2}}: 2^\Theta \to [0, 1.2].$$

So, we divide each mass by 1.1.

Overnormalized $m_{\text{over1,2}}$ of A, B and A∪B respectively gives:

| A | B | A∪B |
|---|---|---|
| $\frac{0.755}{1.1} \simeq$ | $\frac{0.545}{1.1} \simeq$ | $\frac{0.020}{1.1} \simeq$ |
| 0.686 | 0.496 | 0.018 |

Total overmass: 0.686 + 0.496 + 0.018 = 1.2 as required.

The chosen worker for the promotion to manager is A.

## 13. Example of Fusion of OffMasses

We repeat the example of the two workers, A and B, applying for a manager position at a given company, but change the evaluations of the two sources of information (the directors D₁ and D₂).

$$m_{\text{off1}}: 2^\Theta \to [-0.1, 1.2]$$
$$m_{\text{off2}}: 2^\Theta \to [-0.2, 1.3]$$

When we have together overmasses, undermasses, and/or offmasses that get values in various intervals $[a_k, b_k]$, for $a_k < b_k$, we take the *union of these intervals* in order to be able to accommodate all results.

For example, having:

i) $[-0.4, 1]$ and $[0, 1]$, we take $[-0.4, 1]$;

ii) $[0, 1.2]$ and $[-0.1, 1]$, we take $[-0.1, 1.2]$;

iii) $[0, 1.3]$ and $[-0.2, 1]$, we take $[-0.2, 1.3]$;

etc.

**14. Fusion of UnderMasses and OffMasses**

We repeat the example of the two workers, A and B, applying for a manager position at a given company, but change again the evaluations of the two sources of information (the directors $D_1$ and $D_2$).

Since we deal with negative values, we cannot use the conjunctive rule, because by multiplying two negative numbers we get a positive result.

So, we have to use the *Average of Masses Rule*.

Example of UnderMasses:

|  | A | B | A∪B |
|---|---|---|---|
| $m_{under1}$ | −0.2 | 0.7 | 0.3 |
| $m_{under2}$ | 0.4 | −0.1 | 0.5 |

$$m_{under1}: 2^\Theta \to [-0.2, 1]$$
$$m_{under2}: 2^\Theta \to [-0.2, 1]$$
$$m_{under1,2}(A) = \frac{-0.2 + 0.4}{2} = 0.1$$
$$m_{under1,2}(B) = \frac{0.7 + (-0.1)}{2} = 0.3$$
$$m_{under1,2}(A \cup B) = \frac{0.3 + 0.5}{2} = 0.4$$

There is no conflicting mass resulting from this rule. In this case worker B is promoted as manager.

**15. OverMasses in Information Fusion**

Let $\Theta$ be a frame of discernment, and $m: \Theta \to [0, 1]$ a mass (function).

If $\sum_{A \in \Theta} m(A) > 1$, then $m$ is called an *OverMass*.

*Another Example from Information Fusion*

A police investigation is underway, and Ann (A), Brian (B) are two suspects in a crime.

Clues, evidence, and testimonies allow for the establishment of belief masses regarding the guilt of each suspect, based on each available source.

The fusion process allows for the determination of the overall level of belief in the guilt of A or B.

$$m_1, m_2: 2^\Theta \to [0, 1.1], A \cap B = \emptyset, \Theta = A \cup B, \text{FoD} = \{A, B, A \cup B\}.$$

|       | A    | B    | A∪B  | Θ    |                  |
|-------|------|------|------|------|------------------|
| $m_1$ | 0.3  | 0.6  | 0.2  |      |                  |
| $m_2$ | 0.5  | 0.5  | 0.1  |      |                  |
| $m_{12}$ | 0.28 | 0.46 | 0.02 | 0.45 | conjunctive rule |

$$m_{12}(A) = m_1(A)m_2(A) + m_1(A)m_2(A \cup B) + m_2(A)m_1(A \cup B) =$$
$$= 0.3(0.5) + 0.3(0.1) + 0.5(0.2) = 0.28.$$

Similarly for the others. Since

$$\sum_{x \in 2^\Theta} m_1(x) = 0.3 + 0.6 + 0.2 = 1.1.$$

$$\sum_{x \in 2^\Theta} m_2(x) = 0.5 + 0.5 + 0.1 = 1.1.$$

By applying the conjunctive rule, one gets as total sum of masses:

$$[\sum_{x \in 2^\Theta} m_1(x)] \times [\sum_{x \in 2^\Theta} m_2(x)] = 1.1(1.1) = 1.21.$$

Therefore,

$$\sum_{x \in 2^\Theta \cup \emptyset} m_{12}(x) = 0.28 + 0.46 + 0.02 + 0.45 = 1.21.$$

Apply PCR5:

$$\frac{x_A}{0.3} = \frac{y_B}{0.3} = \frac{0.3(0.5)}{0.3 + 0.5} = \frac{0.15}{1.8}$$

$$x_A = \frac{0.3(0.15)}{0.8} = 0.06 \; ; \; y_B = \frac{0.5(0.15)}{0.8} = 0.09$$

$$\frac{x_A}{0.5} = \frac{y_B}{0.6} = \frac{0.30}{1.1} = \frac{0.30}{1.1}$$

$$x_A = \frac{0.5(0.30)}{1.1} = 0.14$$

$$y_B = \frac{0.6(0.30)}{1.1} = 0.16$$

Masses to be added to A: 0.06 + 0.14 = 0.20

*OverNormalized*

Masses to be added to B: 0.09 + 0.16 = 0.25

|                   | A              | B              | A∪B            | Θ     |                                                                                                               |
|-------------------|----------------|----------------|----------------|-------|---------------------------------------------------------------------------------------------------------------|
| $m_{12}$          | 0.28 + 0.20    | 0.46 + 0.25    | 0.02           | 0.45 - 0.45 |                                                                                                               |
| $m_{PCR5}$        | 0.48           | 0.71           | 0.02           | 0.00  | Sum of masses 0.48 + 0.71 + 0.02 + 0.00 = 1.21 OverNormalize the sum to be 1.1 as original. Divide by 1.1 each mass. |
| $m_{PCR5}^{(n)}$  | $\frac{0.48}{1.1}$ = 0.44 | $\frac{0.71}{1.1}$ = 0.64 | $\frac{0.02}{1.1}$ = 0.02 |       | Sum of masses (after overnormalization to 1.1) one has: 0.44 + 0.64 + 0.02 = 1.1                              |

$$Bel(A) = 0.44, Pl(A) = 0.44 + 0.02 = 0.46$$
$$Bel(B) = 0.64, Pl(A) = 0.64 + 0.02 = 0.66$$
$$Bel(A \cup B) = 0.44 + 0.64 + 0.02 = 1.1$$
$$Pl(A \cup B) = 1.1$$

**16. Application Scenario: Sense-Making in Wildfire Evacuation**

In a rapidly evolving wildfire scenario, an Emergency Operations Center (EOC) must synthesize conflicting data from heterogeneous sources, including satellite imagery, IoT ground sensors, and social media. The following scenarios demonstrate how Over/Under/Off Masses provide a superior framework for sense-making compared to classical normalization.

*A. OverMasses for Redundant and Overlapping Evidence*

In emergency response, an OverMass occurs when the sum of belief masses exceeds 1 ($\sum m(A) > 1$). This is common in "hard/soft" fusion where multiple sources report the same event with high confidence, leading to overlapping evidence.

- *The Scenario:* A satellite detects a "high-probability" heat signature ($m$ = 0.9), a ground sensor confirms a "certain" fire ($m$ = 1.0), and social media provides "high-confidence" verification ($m$ = 0.8).
- *Sense-Making:* Instead of classical normalization, which would dampen the evidence to fit the [0, 1] interval, the EOC employs *OverNormalization*. This maintains a sum of masses > 1, signaling "Critical Priority". The intensity of the consensus triggers the immediate deployment of aerial tankers, as the membership degree reflects "over-performance" of evidence similar to overtime work hours.

*B. UnderMasses for Indeterminacy and Coverage Gaps*

An UnderMass occurs when the sum of belief masses is less than 1 ($\sum m(A) < 1$), representing "missing mass" or incomplete information.

- *The Scenario:* In a remote canyon where smoke blocks satellite views and IoT sensors are absent, a single drone feed provides a weak indication of fire spread ($m$ = 0.3).
- *Sense-Making:* The system identifies an UnderMass ($\sum m$ = 0.3, with 0.7 unknown). This triggers a "Search and Reconnaissance" protocol rather than an "Evacuation" protocol. By managing this indeterminacy, the system flags the area for "High Uncertainty," prompting the dispatch of scout teams to fill the gap rather than acting on insufficient data.

*C. OffMasses for Falsification and Security Monitoring*

An OffMass allows for negative mass values ($\sum m(A) < 0$), representing counter-evidence or the active "falsity" of a claim.

- *The Scenario:* A compromised or malfunctioning sensor reports a fire in a "safe zone" ($m$ = 0.8), but thermal cameras and physical observers explicitly confirm the absence of fire ($m$ = -0.8$).
- *Sense-Making:* Classical mass functions struggle to represent "negative" belief. By utilizing OffMasses, the system mathematically "cancels out" the false alarm. This allows the EOC to discount erroneous reports and maintain explainability; the system can justify ignoring the sensor by pointing to the heavy negative weight (damage to credibility) provided by the counter-evidence.

**17. Future Research**

- Fusion of UnderMasses on same intervals
- Fusion of UnderMasses on different intervals
- Fusion of OfMasses on same intervals
- Fusion of OfMasses on different interval
- Fusion of OverfMasses with UnderMass
- Fusion of OverfMasses with OffMass
- Fusion of UnderfMasses with OffMass
- Justification and interpretation of OverMass, UnderMass, and OffMass within the framework of belief functions.

**18. Conclusion**

For the first time we have extended the uncertain (fuzzy, neutrosophic, plithogenic, etc.) OverSet, UnderSet, and OffSet to respectively to OverMass, UnderMass, and OffMass that may have applications in the fusion of information.